\documentclass[11pt,a4paper]{amsart}
\usepackage{amsfonts, amssymb, amsmath, amscd, latexsym}
\usepackage{psfrag, eucal, enumerate, latexsym}
\usepackage{amsthm}
\usepackage{xypic}
\usepackage{enumerate}
\usepackage{psfrag}
\usepackage{graphicx}
\usepackage{amsfonts, amssymb, amsmath, amscd, latexsym}
\usepackage{psfrag, eucal, enumerate, latexsym}
\usepackage{amsthm}
\usepackage{xypic}
\usepackage{enumerate}
\usepackage{psfrag}
\usepackage{graphicx}

\newtheorem{theorem}{Theorem}[section]

\theoremstyle{remark}
\newtheorem{remark}[theorem]{Remark}
\numberwithin{equation}{section}
\theoremstyle{corollary}

\theoremstyle{proposition}
\newtheorem{proposition}[theorem]{Proposition}
\theoremstyle{definition}
\newtheorem{example}[theorem]{Example}

\newcommand{\NN}{\mathbb N}
\newcommand{\PP}{{\mathbb P}}

\newcommand{\ZZ}{{\mathbb Z}}

\newcommand{\OO}{\mathcal{O}}

\DeclareMathOperator{\HH}{H}

\DeclareMathOperator{\rk}{rk}

\DeclareMathOperator{\norm}{norm}
\DeclareMathOperator{\coker}{Coker}

\DeclareMathOperator{\Pic}{Pic}

\title{Semistability of certain bundles on a quintic Calabi-Yau threefold.}

\author{Maria Chiara Brambilla}
\address{ Dipartimento di Matematica e Applicazioni per l'Architettura,\\
Universit\`a di Firenze, piazza Ghiberti, 27\\
50122, Florence, Italy} \email{{brambilla@math.unifi.it}}

\begin{document}
%
%

\begin{abstract}
In the paper ``Chirality change in string theory'',
by Douglas and Zhou, the authors give a list of bundles on a
quintic Calabi-Yau threefold. Here we prove the semistability of
most of these bundles. This provides examples of string theory
compactifications which have a different number of generations and
can be connected.
\end{abstract}

\maketitle
\section{Introduction}
In \cite{fisici} Douglas and Zhou study string theory
compactification and illustrate the chirality change with
different examples.

In particular, in \cite[Section $3$]{fisici} they consider
heterotic string theory with gauge group $E_8 \times E_8$
compactified on a simply-connected compact Calabi-Yau manifold
$M$. In order to show that there exist compactifications on the
same Calabi-Yau with different number of generations  which can be
connected, they need to find examples of semistable holomorphic
vector bundles on the Calabi-Yau, whose Chern classes differ only
in $c_3$.

In \cite[Section $3.3$]{fisici} Douglas and Zhou provide a list of bundles $V$ on a
quintic Calabi-Yau $M\subset\PP^4$, which satisfy the following conditions:
$$c_1(V)=0,\quad c_2(V)=c_2(TM),\quad c_3 \textrm{ arbitrary}.$$
We recall this list in Table \ref{table} above.

Furthermore in order to get supersymmetric vacua it is necessary
to require the semistability of these bundle. In \cite[Appendix
A]{fisici} the authors check one interesting example ($V_8$ in
Table \ref{table}) proving that it is stable against subsheaves
that have a similar monad description.

We recall that a holomorphic vector bundle $V$ on a projective manifold
$X$ with $\Pic(X)\cong\ZZ$ is called {\em stable} if for any
coherent subsheaf $S$ of $V$ with $0<\rk S<\rk V$ we have
$\mu(S)<\mu(V)$, where $\mu=\frac{c_1}{\rk}$, and {\em semistable} if
for any coherent subsheaf $S$ we have $\mu(S)\leq\mu(V)$.
This notion of semistability is also called slope-semistability.

Here we complete the proof of the semistability for most of the
bundles in the list. In particular in Proposition \ref{mia} we
prove  the semistability of the bundles with rank $4$ on a generic
smooth quintic in $\PP^4$. The proof is based on standard
computations and on Flenner's theorem.

In Proposition \ref{tre} we prove that all the sheaves with rank
$3$ in the table, when restricted to a generic smooth quintic
hypersurface in $\PP^4$, are stable bundles. In this case we do
computation directly on the threefold, in order to have locally
free sheaves.

Finally, in Proposition \ref{due} we prove the stability for some
of the bundles with higher rank in the table. In this case we can
prove the stability of a generic bundle with given resolution on a
generic smooth quintic in $\PP^4$. We prove this result by
restricting to a generic plane and using the Dr\'ezet-Le Potier
criterion for the existence of a stable bundle on $\PP^2$. From
this argument we can deduce the stability of our bundles but only
when the resolution is generic in $\PP^2$ (that is only for
bundles $V_{10}$, $V_{12}$, $V_{14}$, $V_{16}$ in Table
\ref{table}).

\section{Preliminaries}
Here we collect some useful results on vector bundles on
projective varieties without giving the proofs. For more details
see e.g.\ \cite{HuLe}.

Let $\PP^n$ denote the complex projective space of dimension $n$.
Let $X$ be a complex projective manifold with $\Pic(X)\cong\ZZ$. A bundle $E$
on $X$ is called {\em normalized} if $c_1(E)\in\{-r+1,\ldots,-1,0\}$,
i.e.\ if  $-1<\mu\leq 0$. We denote by $E_{\norm}$ the unique twist of
$E$ which is normalized.

The following criterion for stability of bundles is a consequence of
the definition:
\begin{proposition}
\label{hop}
Let $V$ be a vector bundle on a projective manifold $X$  with $\Pic(X)\cong\ZZ$.
If $\HH^0(X,(\wedge^qV)_{\norm})=0$ for any $1\leq q\leq \rk(V)-1$, then $V$ is stable.
\end{proposition}

\begin{remark}\label{rem}
Any exact sequence of vector bundles
$$0\to A\to B\to C\to 0$$  induces the following exact sequence for any $q\geq1$
$$0\to S^qA\to S^{q-1}A\otimes B\to\ldots\to A\otimes\wedge^{q-1}B\to\wedge^qB\to\wedge^qC\to0$$
\end{remark}

We state Flenner's theorem in the particular case of hypersurfaces in $\PP^n$:
\begin{theorem}[Flenner]
Assume
$$\binom{d+n}{d}-d-1>d\max\left\{{\frac{r^2-1}{4},1}\right\}.$$ If
$E$ is a semistable sheaf of rank $r$ on $\PP^n$, then the
restriction $E|_X$ on a generic smooth hypersurface $X$ of degree
$d$ in $\PP^n$ is semistable.
\end{theorem}

The following criterion is a particular case (for $c_1=0$) of the
Dr\'ezet-Le Potier theorem (see \cite[Theorem C]{DrezetLePotier}): 
\begin{theorem}[Dr\'ezet-Le Potier]
\label{DLP}
Given $r,c\in\ZZ$ such that
$$c\geq r>0,$$
then there
exists a stable bundle on $\PP^2$ with rank $r$ and Chern classes $c_1=0$ and
$c_2=c$.
\end{theorem}

Let us denote by $M(r,c_1,c_2)$ the moduli space of semistable sheaves
on $\PP^2$ of rank $r$ and Chern classes $c_1,c_2$.   It is known
that $M(r,c_1,c_2)$ is irreducible and moreover we have the following useful
result (see \cite{HiLa} and \cite{DionisiMaggesi}).

\begin{proposition}\label{riso}
A generic bundle in the space $M(r,c_1,c_2)$ has resolution either of
the form
$$
0\rightarrow \OO(k-2)^a\oplus\OO(k-1)^b \rightarrow
 \OO(k)^c\rightarrow F\rightarrow 0,
$$
or
$$
0\rightarrow \OO(k-2)^a \rightarrow
\OO(k-1)^b\oplus \OO(k)^c\rightarrow F\rightarrow 0,
$$
for some $k\in\ZZ$, $a,b\geq0$ and $c>0$.
\end{proposition}
Finally we recall the following 
\begin{proposition}\label{locfree}
Let $\phi:E\to F$ be a morphism of vector bundles on a variety of
dimension $N$, with $e=\rk(E)$, $f=\rk(F)$ and $e\leq f$. If
$E^*\otimes F$ is globally generated and $f-e+1>N$, then for a
generic $\phi$ the sheaf \ $\coker(\phi)$ is locally free, i.e.\
is a vector bundle.
\end{proposition}

\section{Results}
In Table \ref{table} the list of sheaves on $\PP^4$ contained  in
Table $1$ of \cite{fisici} is recalled. To every entry $(n_i,m_j)$
of the table we associate a sheaf with the following resolution
$$0\to V\to\oplus_{i=1}^{r+m}\OO_{\PP^4}(n_i)\to\oplus_{j=1}^{m}\OO_{\PP^4}(m_j)\to0.$$
It is easy to check that all these sheaves on $\PP^4$ have Chern
classes $c_1=0$ and $c_2=10$.

\begin{table}{ \label{table}
\begin{tabular}{|c|c|c|c|}
\hline \hline
&rank & $(n_i)$ & $(m_j)$\\
\hline
\hline
$V_1$&$3$& $(22222222)$&$(33334)$\\
\hline
$V_2$&$3$& $(122222)$&$(344)$\\
\hline
$V_3$&$3$& $(112233)$&$(444)$\\
\hline
$V_4$&$3$&$ (11222)$&$(35)$\\
\hline
$V_5$&$3$& $(11133)$&$(45)$\\
\hline
$V_6$&$4$& $(1122222222)$&$(333333)$\\
\hline
$V_7$&$4$& $(11122222)$&$(3334)$\\
\hline
$V_8$&$4$& $(111122)$&$(44)$\\
\hline
$V_9$&$4$& $(11111)$&$(5)$\\
\hline
$V_{10}$&$5$& $(1111122222)$&$(33333)$\\
\hline
$V_{11}$&$5$& $(11111122)$&$(334)$\\
\hline
$V_{12}$&$6$& $(1111111122)$&$(3333)$\\
\hline
$V_{13}$&$6$& $(111111111)$&$(234)$\\
\hline
$V_{14}$&$7$& $(11111111111)$&$(2333)$\\
\hline
$V_{15}$&$7$& $(111111111111)$&$(22224)$\\
\hline
$V_{16}$&$8$& $(11111111111111)$&$(222233)$\\
\hline
\end{tabular}}
\caption{The list of Douglas and Zhou.}
\end{table}

Moreover if the map in the resolution is generic, by Proposition
\ref{locfree} these sheaves are locally free on $\PP^4$ only if
the rank is bigger or equal than $4$. Nevertheless, by restricting
these sheaves to a generic quintic threefold $M$ in $\PP^4$, we
obtain locally free sheaves also in the case of rank $3$.

 We are interested in proving the (semi)stability of the
restriction of these sheaves to a quintic $M$ in $\PP^4$. First of
all we can prove the following result.

\begin{proposition}\label{mia}
Let $V$ be a generic bundle with rank $4$  in Table \ref{table},
i.e.\ with resolution of type $V_k$, for $6\leq k\leq 9$. Then $V$
is semistable on a generic smooth quintic hypersurface in $\PP^4$.
\end{proposition}

Since (semi)stability is invariant up to duality, we will check
the (semi)\-stability of the dual bundles $V^*$.

Let $E$ be the dual of a bundle of rank $4$ in Table \ref{table}.
To check the stability of $E$ we need to show that
$\HH^0(\PP^4,(\wedge^q E)_{\norm})=0$ for any $1\leq q\leq 3$.
Since $c_1(E)=0$, it is obvious that
$(\wedge^qE)_{\norm}=\wedge^qE$ for any $q$.

\medskip
Before giving the proof of the previous proposition, let us
consider in detail the example of a bundle of rank $4$
corresponding to $V_8$.
\begin{example}\label{cac}
Let $E$ be a bundle with the following resolution on $\PP^4$:
\begin{equation}\label{esempio}
0\to\OO(-4)^2\to\OO(-2)^2\oplus\OO(-1)^4\to E\to0.
\end{equation}
In order to apply Proposition \ref{hop}
we need to check the following conditions
$$
\HH^0(\PP^4,E)=0,\quad
\HH^0(\PP^4,\wedge^2 E)=0,\quad
\HH^0(\PP^4,\wedge^3 E)=0.$$

By the cohomology sequence associated to \eqref{esempio}
we immediately get the first vanishing.
Indeed by Remark \ref{rem}
we can compute the following resolution for $\wedge^2 E$:
$$0\to\OO(-8)^3\to\OO(-6)^4\oplus\OO(-5)^8\to\OO(-4)\oplus\OO(-3)^8\oplus
\OO(-2)^6\to \wedge^2 E\to0,$$
and for $\wedge^3 E$:
$$0\to\OO(-12)^4\to\OO(-10)^6\oplus\OO(-9)^{12}\to\OO(-8)^2\oplus\OO(-7)^{16}\oplus\OO(-6)^{12}\to$$
$$\to\OO(-5)^4\oplus\OO(-4)^{12}\oplus
\OO(-3)^4\to\wedge^3 E\to0.$$
From the resolution of $\wedge^2 E$ we get the following two short exact sequences:
$$0\to\OO(-8)^3\to\OO(-6)^4\oplus\OO(-5)^8\to K_0\to0$$
$$0\to K_0\to \OO(-4)\oplus\OO(-3)^8\oplus
\OO(-2)^6\to \wedge^2 E\to0,$$
and since $\HH^1(\PP^4,K_0)=0$, we get $\HH^0(\PP^4,\wedge^2 E)=0$.
Analogously from the resolution of $\wedge^3 E$ we get
$$0\to\OO(-12)^4\to\OO(-10)^6\oplus\OO(-9)^{12}\to K_1\to0$$
$$0\to K_1\to\OO(-8)^2\oplus\OO(-7)^{16}\oplus\OO(-6)^{12}\to K_2\to0$$
$$0\to K_2\to\OO(-5)^4\oplus\OO(-4)^{12}\oplus\OO(-3)^4\to\wedge^3 E\to0,$$
from which we obtain
$\HH^2(\PP^4,K_1)=0$,
$\HH^1(\PP^4,K_2)=0$, and
$\HH^0(\PP^4,\wedge^3 E)=0$.
Hence we get
 $$\HH^0(\PP^4,E)=\HH^0(\PP^4,\wedge^2 E)=\HH^0(\PP^4,\wedge^3 E)=0,$$
which implies that $E$ is stable on $\PP^4$.

Now from Flenner's theorem it follows that if $E$ is a semistable
bundle of rank $4$ on $\PP^4$, then its restriction on a generic
hypersurface of degree $d\geq2$ is semistable. Hence we conclude
that the restriction of $E$ on a generic smooth quintic
hypersurface in $\PP^4$ is semistable.
\end{example}

\begin{em}{Proof of Proposition \ref{mia}.}\end{em}
It is easy to check that the argument used in Example \ref{cac}
holds for all bundles in Table \ref{table} with rank $4$. Hence
all the bundles with rank $4$ are semistable on a generic smooth
quintic. \qed
\medskip

\begin{proposition}\label{tre}
Let $V$ be a generic sheaf of rank $3$, i.e.\ with resolution of
type $V_k$ for $1\leq k\leq5$. Then the restriction of $V$ to a
generic smooth quintic hypersurface in $\PP^4$ is a stable bundle.
\end{proposition}
\begin{proof}
Let $V$ be any sheaf  of rank $3$ of Table \ref{table}. Since $V$
can be not locally free on $\PP^4$, the argument used in Example
\ref{cac} does not hold. Nevertheless, for  a generic quintic $M$
in $\PP^4$,  the restriction $V|_M$ is locally free, by
Proposition \ref{locfree}.

Let $E|_M$ denote the dual of $V|_M$. For example in the case
$V_5$ we have the following resolution
$$0\to\OO(-5)|_M\oplus\OO(-4)|_M\to(\OO(-3)|_M)^2\oplus(\OO(-1)|_M)^3\to E|_M\to0$$
and we want to apply Proposition \ref{hop} to $E|_M$.

First of all, computing the cohomology of $\OO(-k)|_M$ from the
following exact sequence
$$0\to \OO_{\PP^4}(-5)\to \OO_{\PP^4}\to\OO_M\to 0$$
we easily obtain $\HH^0(M,E|_M)=0$.
On the other hand, by Remark \ref{rem} we can compute the
resolution of $\wedge^2(E|_M)$ and we get
$\HH^0(M,\wedge^2(E|_M))=0$. Analogously we can check that
$$\HH^0(M,E|_M)=\HH^0(M,\wedge^2(E|_M))=0$$
for all cases $V_1,\ldots,V_5$.
 Hence by Proposition \ref{hop} the
restriction $E|_M$ is stable on $M$. We conclude that every sheaf
with rank $3$ of Table \ref{table} restricted to a generic smooth
quintic in $\PP^4$ is a stable bundle.
\end{proof}

Furthermore we can prove that some of the bundles with higher rank
in the table are semistable.
\begin{proposition}\label{due}
Let $V$ be a generic bundle with resolution of the form $V_{10}$,
$V_{12}$, $V_{14}$, or $V_{16}$ in Table \ref{table}. Then $V$ is
semistable on a generic smooth quintic in $\PP^4$.
\end{proposition}
\begin{proof}
Let us show first that $E=V^*$ is stable on $\PP^4$. In order to
do this, it suffices to prove that the restriction of $E$ on a
plane is stable.

If  $V$ is a generic bundle with resolution of the form $V_{10}$,
$V_{12}$, $V_{14}$, or $V_{16}$, then the restriction $E|_{\Pi}$
on a generic plane $\Pi$ has a resolution of the form
$$0\to\OO(-3)^a\oplus\OO(-2)^b\xrightarrow{\phi}\OO(-1)^c\to E|_{\Pi}\to0$$
or
$$0\to\OO(-3)^a\xrightarrow{\phi}\OO(-2)^b\oplus\OO(-1)^c\to E|_{\Pi}\to0,$$
where $a,b,c\in\NN$ are given and $\phi$ is a generic map.

On the other hand, we know from Theorem \ref{DLP} that there
exists a stable bundle on $\PP^2$ with $c_1=0$, $c_2=10$ and
$r\leq 10$. Hence the space $M(r,0,10)$ is not empty for $5\leq
r\leq8$ and irreducible. Hence, by Proposition \ref{riso}, it
follows that a generic bundle in $M(5,0,10)$ (or $M(6,0,10)$,
$M(7,0,10)$, $M(8,0,10)$ respectively) has resolution of the form
$V_{10}$ (or $V_{12}$, $V_{14}$, $V_{16}$ respectively). Therefore
the corresponding bundle $E|_{\Pi}$ is stable on the plane and $E$
is stable on $\PP^4$.

Finally from Flenner's theorem it follows that if $E$ is a
semistable bundle of rank $5\leq r\leq 8$ on $\PP^4$, then its
restriction to a generic hypersurface of degree $5$ is semistable.
Hence we conclude that the restriction of $E$ on a generic smooth
quintic hypersurface in $\PP^4$ is semistable.
\end{proof}
\medskip
\begin{bf}{Acknowledgments}\end{bf} {The author thanks Adrian Langer for
suggesting this problem and for useful discussions, and Giorgio
Ottaviani for helpful comments. The author was partially supported
by Italian MIUR funds and by a CNR-NATO Advanced Fellowship.}


\begin{thebibliography}{1}

\bibitem{DionisiMaggesi}
Carla Dionisi and Marco Maggesi, \emph{Minimal resolution of
general stable
  rank-2 vector bundles on {${\PP}^2$}}, \emph{Boll. Unione Mat. Ital. Sez. B Artic.
  Ric. Mat.} (8) \textbf{6} (2003), no.~1, 151--160,
  math.AG/0002169.

\bibitem{fisici}
Michael~R. Douglas and Chengang Zhou, \emph{Chirality change in
string theory}, \emph{ J. High Energy Phys.} (2004), no.~6, 014,
42 pp. (electronic).

\bibitem{DrezetLePotier}
Jean-Marc Dr{\'e}zet and Joseph Le~Potier, \emph{Fibr\'es stables
et fibr\'es
  exceptionnels sur {${\PP}^2$}}, \emph{Ann. Sci. \'Ecole Norm. Sup.} (4) \textbf{18}
  (1985), no.~2, 193--243.

\bibitem{HiLa}
Andr{\'e} Hirschowitz and Yves Laszlo, \emph{Fibr\'es g\'en\'eriques sur le
  plan projectif}, \emph{Math. Ann.} \textbf{297} (1993), no.~1, 85--102.

\bibitem{HuLe}
Daniel Huybrechts and Manfred Lehn, \emph{The geometry of moduli
spaces of sheaves}, {Aspects of Mathematics, \textbf{31}}, Vieweg,
{Braunschweig}, Germany 1997.





\end{thebibliography}
\end{document}